\documentclass[11pt,reqno]{amsart}

\usepackage{eucal,color,graphicx,mathrsfs,tikz}
\usepackage{hyperref}
\usepackage{amsthm}

\usepackage{mathtools}

\renewcommand{\Omega}{\om}

\newcommand{\om}{\varOmega}

\newtheorem{theorem}{Theorem}[section]

\newtheorem{corollary}[theorem]{Corollary}
\newtheorem{proposition}[theorem]{Proposition}

\theoremstyle{definition}
\newtheorem{definition}[theorem]{Definition}
\newtheorem{example}[theorem]{Example}

\theoremstyle{remark}

\numberwithin{equation}{section}

\begin{document}

\title{Zheghalkin-Boolean Calculus}

\author[S. Nagaraj]{Sriram ~Nagaraj}
\email{sriram.nagaraj.atl@gmail.com}
 
\thanks{Corresponding author: Sriram Nagaraj (sriram.nagaraj.atl@gmail.com), with the Federal Reserve Bank of Atlanta. The views expressed here are the author's and not necessarily those of the Federal Reserve Bank of Atlanta or the Federal Reserve System.}

\begin{abstract}
Boolean calculus has been studied extensively in the past in the context of switching circuits, error-correcting codes etc. This work generalizes several approaches to defining a differential calculus for Boolean functions. A unified theory of Boolean calculus, complete with $k$-forms and integration, is presented through the use of Zhegalkin algebras (i.e., algebraic normal forms), culminating in a Stokes-like theorem for Boolean functions. 
\end{abstract}
\maketitle
\section{Introduction}
The ``calculus" of Boolean functions, i.e., functions with $n$-binary inputs and a binary output, has been studied, in one form or another, since the time of Reed \cite{reed}, Muller \cite{muller}, Huffman \cite{huffman} and other early pioneers for the design and testing of switching circuits, error-correcting codes etc. Different approaches to a theory of calculus for Boolean functions have been proposed, each motivated by specific applications (such as circuit optimization, equivalent reduced sum-of-products/product-of-sums forms etc.). In this paper, the main motivation is to develop a sufficiently general foundation that allows us to state and prove a Stokes-like theorem, thus justifying the term ``calculus".
\subsection{Prior Work}
The work in this paper is motivated by the construction \cite{zheg} of Zhegalkin (or Gegalkine, as the name is sometimes spelled) and are one of the possible representations of the operations of a general Boolean algebra, and correspond to the algebraic normal form (ANF) of a Boolean function. These monomial algebras were introduced by I. I. Zhegalkin in 1927.
In the course of several decades since the original motivational work, various authors have proposed several related constructions of a differential calculus and associated theory of differential equations for Boolean functions \cite{StPos1,StPosBook1,StPosBook2,sefarati95,thayse}. Other developments include a theory of integral calculus for Boolean functions \cite{tucker}, as well as general theories of differential operators/algebras \cite{aragona,encinas,ruiz,kuhnrich} in the Boolean context. Boolean calculus with a linear algebraic flavor was presented in \cite{yanu,mizraji}, and a matrix theory approach, including indefinite integrals and primitive functions, was provided in \cite{boolmat}. The work of \cite{encinas} also uses polynomials (algebras), in the Boolean context.
\subsection{Contributions}
This paper addresses the need for a unified theory of Boolean calculus. Indeed, the general version of Stokes theorem \cite{spivak} can be viewed as a unified statement of the theory of (classical) calculus. The aim of this paper is to present a theory of calculus of Boolean functions which is sufficiently general to be able to prove a Stokes theorem like result in the Boolean case. Since Boolean algebras are, in some sense, devoid of non-trivial derivations \cite{rudeanu}, a technical alternative is provided which suffices to allow our theory to go through.

\subsubsection{Outline of Paper}
After fixing some notation, section \ref{sec:1} is devoted to defining the Zhegalkin algebras of degree $n$ and studying the (lack of) derivations of these algebras. In lieu of a space of derivations (i.e., tangent vectors), the space of secants is defined as a proxy for performing differentiation. The corresponding dual space of $1$-forms is also defined. Section \ref{sec:2} defines the space of $k$-forms on the dual secant space and studies the interaction of the wedge product with the proxy-derivative map. Section \ref{sec:3} defines the integration theory of $k$-forms and sets the stage for the Stokes-Zhegalkin theorem.

\section{Zhegalkin Algebras and the Space of Secants}\label{sec:1}
\subsection{Notation}
We now fix the main notation to be used in the remainder of this article. Other specific notation will be introduced when requried. All rings considered in this article will be commutative. A \emph{Boolean ring} is a ring $(R,+,\cdot)$ equipped with two binary operations $+,\cdot : R\times R \rightarrow R$ such that $a\cdot a = a$ for all $a$ in $R$. Note that a Boolean ring is always commutative and $a+a=0$ for all $a  \in R$ as well. We remind the reader that any Boolean ring can alternatively be described as a Boolean algebra, i.e., a  complemented distributive lattice $(B, \wedge,\vee,\neg)$ consisting of a set $B$ with two binary operations, namely conjunction and disjunction, denoted by $\wedge,\vee: B\times B\rightarrow B$, and a unary operation, namely negation, denoted by $\neg: B\rightarrow B$. Indeed, given a Boolean ring $(R, +,\cdot)$, we can define, for any $a,b\in R$:\begin{align} 
a\wedge b &:=  a\cdot b \\ 
a\vee b &:=  a+b+a\cdot b \\
\neg a &:= 1+a \\
\end{align}that turns $(R, \wedge,\vee,\neg)$ into a Boolean algebra. One can analogously turn a Boolean algebra into a Boolean ring.

For any set $X$, we define $I_X: X \rightarrow X$ to be the identity map. Given a ring $R$, and indeterminate generators $x_1,\ldots,x_n$, we denote by $$\text{span}_{R}\{x_1,\ldots,x_n\}$$ the free $R$ module generated by the elements $x_1,\ldots,x_n$. Given an ideal $I \subset R$, $R/I$ will denote the quotient ring of $R$ by $I$. The direct sum of rings (resp. modules) $R_1,R_2$ will be denoted by $R_1 \oplus R_2$. Given $R$ modules $M_1,M_2$, their $R$-tensor product will be denoted by $M_1 \otimes_{R} M_2$. The isomorphic relationship between two objects (rings, modules) will be denoted by $\cong$, so that $X_1 \cong X_2$ means $X_1$ and $X_2$ are isomorphic objects.
Given a set of indeterminates $\{x_1,\ldots,x_n\}$ and a ring $R$, we denote $R[x_1,\ldots,x_n]$ to be the polynomial ring over the indeterminates $x_i, \, i=1,\ldots,n$ with coefficients in $R$. We refer the reader to \cite{df,matsumura} for further algebraic details. 
The ring of integers will be denoted by $\mathbb{Z}$ and the field of integers modulo $2$ will be denoted by $\mathbb{F}_2$, i.e., $$\mathbb{Z}/2\mathbb{Z} := \mathbb{F}_2 = \{0,1\},$$ and $\mathbb{F}_{2}^{n}$ is the $n$-dimensional $\mathbb{F}_2$ vector space obtained by taking the $n$-fold direct sum of $\mathbb{F}_2$, in other words, $$\mathbb{F}_{2}^{n} := \bigoplus_{i=1}^{n} \mathbb{F}_{2}.$$ In the present case, $\mathbb{F}_{2}^{n}$ is a Boolean ring with componentwise addition and multiplication. A \emph{Boolean} function is any map $$f:\mathbb{F}_{2}^{n} \rightarrow \mathbb{F}_{2}.$$Note that the set of all Boolean functions has cardinality $2^{2^{n}}$.
A fundamental result (see for e.g. \cite{rost}) that motivates our subsequent discussion is as follows.
\begin{theorem}Any Boolean function $f$ can be represented uniquely as a sum of monomials, $$f(x_1,\ldots,x_n) = \sum_{0 \leq i_1,\ldots,i_k \leq n}a_{i_{1}}\ldots a_{i_{k}}x_{1}^{m}\ldots x_{n}^{m} \in \mathbb{F}_{2}[x_1,\ldots,x_n],$$ with $m \in \{0,1\}$, $x_{i}^{0}=1, x_{i}^{1}=x_{i}$ and coefficients $a_{i_{j}}\in \mathbb{F}_{2}$ for $i,j=1,\ldots,n$.
\end{theorem}
In short, the calculus we shall build rests on the above result on monomial representations of Boolean functions, and we devote the following sections of the paper to a detailed study of these representations.
\subsection{The Zhegalkin Algebra $\mathfrak{Z}_n$}

We reviewed the relationship between the space of Boolean functions and their monomial representation in the end of previous subsection. We shall formalize these monomial representations by defining and studying the \emph{Zhegalkin algebra} of degree $n$ shortly.

Let $\mathbb{F}_2 \langle x \rangle := \mathbb{F}_2[x]/(x^2+x)$ where $(x^2+x)\subset \mathbb{F}_2[x]$ is the ideal generated by the $\mathbb{F}_2$ polynomial $x^2+x$. Notice that $\mathbb{F}_2 \langle x \rangle$ is a two dimensional torsion module over $\mathbb{F}_2$ generated by $\{1,x\}$ subject to the relation $x^2=x$. Thus, $\mathbb{F}_2 \langle x \rangle$ is the prototypical Boolean ring, and is isomorphic (as an $\mathbb{F}_2$ module) to $\mathbb{F}_2\oplus\mathbb{F}_2$. We can now define the \emph{Zhegalkin algebra} of degree $n$.
\definition{Given indeterminates $x_1,\ldots,x_n$, the \emph{Zhegalkin algebra of degree $n$} is defined as:\begin{equation}\mathfrak{Z}_n = \bigotimes_{i=1}^{n}\mathbb{F}_2 \langle x_i \rangle = \mathbb{F}_2 \langle x_1 \rangle \otimes_{\mathbb{F}_2}\mathbb{F}_2 \langle x_2 \rangle \ldots \otimes_{\mathbb{F}_2}  \mathbb{F}_2 \langle x_n \rangle.\end{equation}}

We notice that $\mathfrak{Z}_n$ is also a Boolean ring. In fact, by construction, it is a $2^n$ dimensional $\mathbb{F}_2$ vector space and inherits its vector space structure from the constituent $\mathbb{F}_2 \langle x_i \rangle$. Moreover, we see that the map $i:\mathbb{F}_2\hookrightarrow \mathfrak{Z}_n$ with $$i(a) = a(\underbrace{1 \otimes 1 \otimes \ldots \otimes 1}_{n\text{ factors}}),$$ for $a \in \mathbb{F}_2$ is an embedding of $\mathbb{F}_2$ into $\mathfrak{Z}_n$. We now have the following proposition.

\begin{proposition}Let $S=\text{span}_{\mathbb{F}_2}\{x_{1}^{\alpha_1}\ldots x_{n}^{\alpha_1}\}$ where $\alpha_i \in \mathbb{F}_2, \, i=1,\ldots,n$ and $x_{i}^{1}=x_{i}$ while $x_{i}^0=1$. Thus, $S$ is the $\mathbb{F}_2$ span of monomials. As modules, $\mathfrak{Z}_n$ and $S$ are isomorphic.
\begin{proof}
Notice first that $\mathbb{F}_2 \langle x_i \rangle = \text{span}_{\mathbb{F}_2}\{x_{i}^{\alpha_i} :\alpha_i \in \mathbb{F}_2\}$. Thus, every $x\in \mathfrak{Z}_n$ can be written as a sum of terms of the form $$x=ax_{1}^{\alpha_1}\otimes x_{2}^{\alpha_2}\otimes \ldots \otimes x_{n}^{\alpha_n},$$ where $a \in \mathbb{F}_2$. We define $\phi: \mathfrak{Z}_n \rightarrow S$ first on monomials and extend by linearity. Thus $\phi(x) = ax_{1}^{\alpha_1}x_{2}^{\alpha_2}\ldots x_{n}^{\alpha_n} \in S$. Clearly $\phi$ is surjective by definition. If $\phi(x)=0$, then either $a=0$ or $x_i=0$ for some $i$. In either case, $x=0$. Thus, $\phi$ is injective.

Before extending by linearity, we verify $\phi(xy)=\phi(x)\phi(y)$. This is a routine calculation. If $y=bx_{1}^{\beta_1}\otimes x_{2}^{\beta_2}\otimes \ldots \otimes x_{n}^{\beta_n}$ with $\beta_i\in \mathbb{F}_2, \, i=1,\ldots,n$, then $$xy = ab (x_{1}^{\alpha_1} x_{1}^{\beta_1}\otimes x_{2}^{\alpha_2} x_{2}^{\beta_2}\otimes \ldots x_{n}^{\alpha_n} x_{n}^{\beta_n}),$$ and hence, $$\phi(xy) = ab (x_{1}^{\alpha_1} x_{1}^{\beta_1}x_{2}^{\alpha_2} x_{2}^{\beta_2} \ldots x_{n}^{\alpha_n}x_{n}^{\beta_n}) =  (a \, x_{1}^{\alpha_1} x_{2}^{\alpha_2} \ldots x_{n}^{\alpha_n})(b\, x_{1}^{\beta_1} x_{2}^{\beta_2}\ldots  x_{n}^{\beta_n}) = \phi(x)\phi(y).$$

Finally, we extend $\phi$ to all of $\mathfrak{Z}_n$ by linearity and the result follows.
\end{proof}
\end{proposition}

Thus, we see that $\mathfrak{Z}_n$ can be identified with the free span of the monomials $\text{span}_{\mathbb{F}_2}\{x_{1}^{\alpha_1}\ldots x_{n}^{\alpha_1}\}$ where $\alpha_i \in \mathbb{F}_2, \, i=1,\ldots,n$ subject to the relations $x_{i}^{2}=x_{i}$ (so that $x_i+x_i=0$ as well). For the remainder of this paper, we will view the algebra $\mathfrak{Z}_n$ as this span of monomials. As an example, $\mathfrak{Z}_2 = \text{span}_{\mathbb{F}_2}\{1\otimes 1, x_1 \otimes 1, 1 \otimes x_2, x_1 \otimes x_2\}$. Finally, we see that, as mentioned earlier, any Boolean function of $n$ variables can be represented by elements of $\mathfrak{Z}_n$.

Having defined the Zhegalkin algebras $\mathfrak{Z}_n$, we proceed to studying certain canonical linear functionals defined on $\mathfrak{Z}_n$ that will be used to build the Zhegalkin-Boolean calculus.

\subsection{Canonical Secants of $\mathfrak{Z}_n$}

In classical calculus, the notion of the \emph{tangent space} plays a central role. Defining the tangent space at a point on a smooth manifold, for example, can be done through several equivalent ways. One common way is through the observation that tangent vectors are derivations on the ring of (germs of) smooth function defined at the point in question. 

In the present case of Zhegalkin algebras, with $\mathfrak{Z}_n$ being the ring of functions on which one looks for derivations, this approach does not work well. Indeed, any reasonable definition of derivations on a Boolean ring lead to trivialities. We refer the reader to \cite{rudeanu} for further details and analysis. However, there is also a simple ``sanity check" that leads us to the same conclusion.

We consider the case of $\mathfrak{Z}_1 = \mathbb{F}_2 \langle x_1 \rangle$. A classical way of defining the algebraic tangent space is through the module of \emph{K\"ahler differentials} (see \cite{matsumura} for details). Briefly, the module of K\"ahler differentials of the ring $R[x_1,\ldots,x_n]/(f_1,\ldots,f_m)$ can be identified with the free $R[x_1,\ldots,x_n]$-module spanned by differentials $dx_1,\ldots,dx_n$ modulo the ideal $(df_1,\ldots,df_m)$. In the case of $\mathfrak{Z}_1 = \mathbb{F}_2[x_1]/(x_{1}^{2}+x_1)$, we see that the module of K\"ahler differentials of $\mathfrak{Z}_1 =0$ since the differential ideal of $(x_{1}^{2}+x_1)$ is the entire algebra spanned by $dx_1$. Thus, the classical cotangent (and hence tangent) space of $\mathfrak{Z}_1=0$. Tensoring this observation $n$ times yields the fact that the tangent space of $\mathfrak{Z}_n=0$.

In order to continue with our theory, we proceed as follows. Consider again the case of $\mathfrak{Z}_1$ first. We start with the space of linear functionals on $\mathfrak{Z}_1$ into $\mathbb{F}_2 $ and extend the range to all of $\mathfrak{Z}_1$ using the embedding $i: \mathbb{F}_2 \hookrightarrow \mathfrak{Z}_1$. The space of linear functionals on $\mathfrak{Z}_1$ contains in particular, the mapping $e_1$ that takes $1 \mapsto 0$ and $x_1 \mapsto 1$. It is then extended by linearity. Finally, composing $e_1$ with the embedding $i: \mathbb{F}_2 \hookrightarrow \mathfrak{Z}_n$, we obtain our ``derivative surrogate" $e_{1}^{'} = i \circ e_1$. We refer to the distinguished linear map $e_{1}^{'}$ as a \emph{secant} in the space of all endomorphisms of $\mathfrak{Z}_1$. 

The same construction performed on $\mathbb{F}_2 \langle x_i \rangle, \, i=2,\ldots,n$ gives us a total of $n$ secants $e_{i}^{'}, \, i=1,\ldots,n$ defined on $\mathbb{F}_2 \langle x_i \rangle$ respectively.  Each $e_{i}^{'}$ can be further extended to a secant on $\mathfrak{Z}_n$ as follows:
\begin{equation}\partial_i = \underbrace{1 \otimes 1 \otimes \ldots \otimes \overset{\substack{\text{i-th factor}\\\uparrow}}{e_{i}^{'}} \otimes 1 \ldots \otimes 1}_{n\text{ factors}}\end{equation}
Thus the maps $\partial_i: \mathfrak{Z}_n \rightarrow \mathfrak{Z}_n$ are linear and correspond to ``partial differentiation" with respect to $x_i$ for $i_1,\ldots,n$ with the caveat that they are not true derivations on the algebra $\mathfrak{Z}_n$. Note that the, since $\mathfrak{Z}_n$ is comprised of monomials, the image of $\partial_i$ is the sub-algebra of $\mathfrak{Z}_n$ which contains no terms with an $x_i$ factor.

\subsection{The Secant Space $\mathcal{S}(\mathfrak{Z}_n)$}
The collection of secants $\partial_i:\mathfrak{Z}_n \rightarrow \mathfrak{Z}_n$ are linearly independent in the space of endomorphisms of $\mathfrak{Z}_n$ by construction. Indeed, since the image of $\partial_i$ and $\partial_j$ for distinct $i,j$ contain no terms with $x_i$ and $x_j$ respectively, there is no non-trivial linear relationship between $\partial_i$ and $\partial_j$. We can thus consider the $n$-dimensional algebra $\mathcal{S}(\mathfrak{Z}_n)=\text{span}_{\mathfrak{Z}_n}\{\partial_1,\ldots,\partial_n\}$ defined by the secants $\partial_i$ over the algebra $\mathfrak{Z}_n$. We define\begin{definition}The $n$-dimensional $\mathfrak{Z}_n$-algebra $\mathcal{S}(\mathfrak{Z}_n) = \text{span}_{\mathfrak{Z}_n}\{\partial_1,\ldots,\partial_n\}$ is defined to be the secant space generated by the secants $\partial_i:\mathfrak{Z}_n \rightarrow \mathfrak{Z}_n$.\end{definition}
Given any $\phi \in \mathcal{S}(\mathfrak{Z}_n)$, we can express $\phi$ as $\phi = \sum_{i=1}^{n}f_{i}\partial_i$ with $f_{i} \in \mathfrak{Z}_n$. In addition, the map $\phi$ acts on any $g\in \mathfrak{Z}_n$ as $\phi(g)=\sum_{i=1}^{n}f_{i}\partial_i(g)\in \mathfrak{Z}_n$. Moreover, any $g\in \mathfrak{Z}_n$ acts on $\phi$ as $g(\phi) = (g\phi) \in \mathcal{S}(\mathfrak{Z}_n)$. It is easy to see these actions are linear, i.e., for $f,g \in \mathfrak{Z}_n$ and $\phi,\psi \in \mathcal{S}(\mathfrak{Z}_n)$, we readily observe that $\phi(f+g) = \phi(f)+\phi(g)$, $(\phi+\psi)(f) = \phi(f)+\psi(f)$, and finally that $(f+g)(\phi)=f\phi+g\phi$.
We can also define the (algebraic) dual of $\mathcal{S}(\mathfrak{Z}_n)$ to be the space $\Omega(\mathfrak{Z}_n)$ of 1-forms:
\begin{definition}The space $\Omega(\mathfrak{Z}_n)$ of 1-forms is defined to be the algebraic dual $(\mathcal{S}(\mathfrak{Z}_n))^{*}$ of $\mathcal{S}(\mathfrak{Z}_n)$.\end{definition}
Given that the $\partial_i$ span $\mathcal{S}(\mathfrak{Z}_n)$, there exist the corresponding dual basis $d_i$ that span $\Omega(\mathfrak{Z}_n)$. Indeed, the basis $d_i$ of $\Omega(\mathfrak{Z}_n)$ is characterized by the relations $d_i(\partial_j)=\delta_{ij}$ where $\delta_{ij}$ is the Kronecker symbol: $\delta_{ij} = 1$ if $i=j$ and is $0$ otherwise. Finally, we can define the map $d:\mathfrak{Z}_n \rightarrow \Omega(\mathfrak{Z}_n)$ as follows.
\begin{definition}Define the map $d:\mathfrak{Z}_n \rightarrow \Omega(\mathfrak{Z}_n)$ to be such that $df(\phi)=\phi(f)$ for any $f\in \mathfrak{Z}_n$ and $\phi\in\mathcal{S}(\mathfrak{Z}_n)$.\end{definition}
The map $d$ defined above will play a pivotal role in the development of our Boolean calculus. It is clear that $d$ is also linear since $df(\phi+\psi)=(\phi+\psi)(f)=\phi(f)+\psi(f)=df(\phi)+df(\psi)$ for $f \in \mathfrak{Z}_n$ and $\phi,\psi \in \mathcal{S}(\mathfrak{Z}_n).$
We close this section with the following observation. We can apply the map $d$ to the monomial indeterminates $x_i \in \mathfrak{Z}_n$, and obtain the analogs of the ``differentials" $dx_i$ of the ``usual" calculus. Indeed, $dx_i(\partial_j) = \partial_j(x_i) = \delta_{ij}$ and hence, we observe that $dx_i=d_i$. We can now establish the following result.

\begin{proposition}Let $f\in \mathfrak{Z}_n$. Then, $df = \sum_{i=1}^{n}\partial_{i}(f) d_{i}$.
\begin{proof}
Given any $f\in \mathfrak{Z}_n$, we have that $df\in \Omega(\mathfrak{Z}_n)$ and hence, there exist $g_i\in \mathfrak{Z}_n$ such that $df = \sum_{i=1}^{n}g_{i} d_{i}$. Applying $\partial_j$ to both sides of the previous relation, we obtain $df(\partial_j) = \partial_j(f) = \sum_{i=1}^{n}g_{i}d_i(\partial_j) = \sum_{i=1}^{n}g_{i}\delta_{ij}=g_j$. Thus $g_j=\partial_j(f)$, and hence $df = \sum_{i=1}^{n}\partial_{i}(f) d_{i}$.
\end{proof}
\end{proposition} 

\section{The Space of $k$-Forms}\label{sec:2}
This section is devoted to the development of the analogs of $k$-forms in the Boolean setting. We have already seen the space of $1$-forms $\Omega(\mathfrak{Z}_n)$ in the previous section. Since $\Omega(\mathfrak{Z}_n)$ is a $\mathfrak{Z}_n$ module, we can construct the alternating module $\Lambda(\Omega(\mathfrak{Z}_n))$ of $\Omega(\mathfrak{Z}_n)$ \cite{matsumura,df}. This alternating structure $\Lambda(\Omega(\mathfrak{Z}_n))$ comes with the alternating (``wedge") product $\wedge: \Lambda(\Omega(\mathfrak{Z}_n)) \rightarrow \Lambda(\Omega(\mathfrak{Z}_n))$ and an associated grading, i.e., a collection of graded submodules $\Lambda^{k}(\Omega(\mathfrak{Z}_n)), k=0,\ldots,n$ such that $\Lambda(\Omega(\mathfrak{Z}_n))$ is the direct sum of the $\Lambda^{k}(\Omega(\mathfrak{Z}_n))$, i.e., 
\begin{equation}\Lambda(\Omega(\mathfrak{Z}_n))=\bigoplus_{k=0}^{n}\Lambda^{k}(\Omega(\mathfrak{Z}_n)),\end{equation}
and
\begin{equation}\mathfrak{Z}_{n}= \Lambda^{0}(\Omega(\mathfrak{Z}_n)) \subset \Lambda^{1}(\Omega(\mathfrak{Z}_n)) \subset \ldots \subset \Lambda^{n}(\Omega(\mathfrak{Z}_n)).\end{equation}
In addition, the wedge product $\wedge: \Lambda^{k}(\Omega(\mathfrak{Z}_n)) \rightarrow \Lambda^{k+1}(\Omega(\mathfrak{Z}_n)).$ Finally, the dimension of $\Lambda^{k}(\Omega(\mathfrak{Z}_n))$ as a free module is ${n \choose k}$. The elements $\omega$ of $\Lambda^{k}(\Omega(\mathfrak{Z}_n))$ are called $k$-forms and are represented as $\omega = \sum_{i_{1}\ldots i_{k}}f_{i_{1}\ldots i_{k}}d_{i_{1}}\wedge d_{i_{2}} \wedge \ldots d_{i_{k}}$ with $f_{i_{1}\ldots i_{k}} \in \mathfrak{Z}_n$. Note that, in contrast with the usual calculus, the wedge product is commutative due to the base field being $\mathbb{F}_2$: for $\omega,\eta \in \Lambda(\Omega(\mathfrak{Z}_n))$, we have $\omega \wedge \eta = -\eta \wedge \omega = \eta \wedge \omega$.

We now study the relation of $\wedge$ with the $d:\Omega(\mathfrak{Z}_n) \rightarrow \Omega(\mathfrak{Z}_n)$ map.
\begin{definition}Extend the map $d$ defined on $\mathfrak{Z}_n$ to the space of $k$ form $\Lambda^{k}(\Omega(\mathfrak{Z}_n))$ as follows. For $\omega \in \Lambda^{k}(\Omega(\mathfrak{Z}_n))$ with $\omega = \sum_{i_{1}\ldots i_{k}}f_{i_{1}\ldots i_{k}}d_{i_{1}}\wedge d_{i_{2}} \wedge \ldots d_{i_{k}}$ define \begin{equation}d\omega = \sum_{i_{1}\ldots i_{k}}df_{i_{1}\ldots i_{k}}\wedge d_{i_{1}}\wedge d_{i_{2}} \wedge \ldots d_{i_{k}} \in \Lambda^{k+1}(\Omega(\mathfrak{Z}_n)).\end{equation}Thus $d:\Lambda^{k}(\Omega(\mathfrak{Z}_n)) \rightarrow \Lambda^{k+1}(\Omega(\mathfrak{Z}_n))$\end{definition}
Note that while in the ordinary calculus, we have the ``product rule" $d(\omega \wedge \eta) = d\omega \wedge \eta + (-1)^{k}\omega \wedge d\eta $, such a relationship does not exist in the present case due to the fact that the underlying map $d$ is not a derivation. However, it is easily seen that $d^2=0$.

\section{Integration of Forms: Stokes-Zheghalkin Theorem}\label{sec:3}
For $f\in \mathfrak{Z}_n$, define the evaluation $f|_{x_i=b_i}=f(x_1,\ldots,x_i=b_i,x_{i+1},\ldots,x_n) \in \mathfrak{Z}_n$ for $b_i\in \mathbb{F}_2$. Thus $f|_{x_i=b_i}$ evaluates the function $f$ at $x_i=b_i$ holding the other arguments $x_j, j\neq i$ fixed and yields a function devoid of an $x_i$ argument. Likewise, $f|_{x_i=b_i,x_j=b_j}$ indicates that $f$ is evaluated at both $x_i=b_i$ and $x_j=b_j$ thereby yielding a function hving no $x_i,x_j$ arguments. Thus, $f|_{x_1=b_1,x_2=b_2,\ldots,x_n=b_n}$ yields a constant obtained by evaluating $f$ at $x_1=b_1,x_2=b_2,\ldots,x_n=b_n$. We start with the following observation. Compare this with the results of \cite{encinas,thayse,StPosBook1} to see the obvious similarities.
\begin{proposition}\label{prop:I1}Let $f\in \mathfrak{Z}_n$. Then $\partial_i(f) = f|_{x_i=0}+f|_{x_i=1}.$
\begin{proof}
We decompose $f$ as $f=f_1+x_if_2$ where the $f_1,f_2 \in \mathfrak{Z}_n$ do not have any $x_i$ factors. This is always uniquely possible since $f$ is a product of monomial terms. Now, on the one hand, $\partial_i(f) = f_2$. On the other hand, we see that $f|_{x_i=0} = f_1$ and $f|_{x_i=1} = f_1+f_2$, so that $f|_{x_i=0} + f|_{x_i=1} = f_1+f_2+f_1=f_2$. The result follows.
\end{proof}
\end{proposition}
We immediately obtain:
\begin{corollary} Let $f\in \mathfrak{Z}_n$. Then $df=(\sum_{i=1}^{n}f|_{x_i=0} + f|_{x_i=1})d_i$
\begin{proof}
Follows from the fact that $df=\sum_{i=1}^{n}\partial_i(f) d_i$ and proposition \ref{prop:I1}
\end{proof}
\end{corollary}
\subsection{Hamming Cube}
We define the Hamming cube $\mathbb{H}^n = \mathbb{F}^n$, which we view as the vertices of the unit cube in $n$-dimensions. Note that $\mathbb{H}^n$ is comprised of $2^n$ binary vectors, each of which encode the binary representation of integers $0\le k \le 2^{n}-1$ i.e., $$\mathbb{H}^n =\{v_k \in \mathbb{F}^n: \text{entry $j$ of }v_k \text{ is the $j$-th entry of the binary expansion of }k,  0\le k \le 2^{n}-1\}.$$
We define the $(i,j)$ boundary $\mathbb{H}^{n-1}_{i,j}, i \in \mathbb{F}_{2}^{n}, j \in \mathbb{F}_2$ of $\mathbb{H}^n$ as:
\begin{equation}\mathbb{H}^{n-1}_{i,j} = \{v_k \in \mathbb{H}^n: (v_k)_{i}= j\}.\end{equation}
Thus, $\mathbb{H}^{n-1}_{i,j}$ contains all binary vectors in $\mathbb{H}^n$ whose $i$-th component is $j$, and is therefore a ``face" of the Hamming cube.
Finally, we define the boundary $\partial \mathbb{H}^n$ of $\mathbb{H}^n$ as the (disjoint) union of all the faces of the cube, i.e., $(i,j)$ boundaries:
\begin{equation}\partial \mathbb{H}^{n-1} = \bigsqcup_{i=1,\ldots,n, \, j\in \mathbb{F}_2}\mathbb{H}^{n-1}_{i,j}.\end{equation}

\subsection{Integration of forms}
In defining the integral of a general $k$-form, henceforth, we restrict our attention to ``monomial" forms, i.e., those that can be expressed as $\omega = f d_{i_1}\wedge \ldots \wedge d_{i_k}$. Although we will primarily be dealing with $n$ and $(n-1)$-forms in this section, we define the \emph{support} of a monomial $k$-form for $0\le k\le n$. Given an $k$ form, say $\omega = f d_{i_1}\wedge \ldots \wedge d_{i_k}$, we define its support to be the subset of $\mathbb{H}^n$ given by \begin{equation}supp(\omega) = \bigcup_{j=k+1,\, b_j \in \mathbb{F}_2}^{n} \mathbb{H}^{n-1}_{i_{k+1},0}\cap \mathbb{H}^{n-1}_{i_{k+2},0} \cap \ldots \cap \mathbb{H}^{n-1}_{i_{j},b_j}\ldots \cap \mathbb{H}^{n-1}_{i_{n},0}\end{equation}Thus, for an $n$-form, its support is defined to be all of $\mathbb{H}^n$. If $\omega=f d_1 \wedge d_2 \ldots \wedge \hat{d}_k \wedge \ldots d_n$ is an $(n-1)$-form (where $\hat{d}_{k}$ indicates the absence of the quantity $\hat{d}_{i}$), its support is the set consisting of the two faces $\mathbb{H}^{n-1}_{k,0},\mathbb{H}^{n-1}_{k,1}$. Now, let $\omega = f d_{i_1}\wedge \ldots \wedge d_{i_k}$ be a $k$-form. We define the integral of $\omega$ over its support as follows:
\begin{definition}$\int_{supp(\omega)}\omega =\sum_{b_{i_1},b_{i_2},\ldots,b_{i_k} \in \mathbb{F}_2}(x_{i_{1}} x_{i_{2}}\ldots x_{i_{k}} f)|_{x_{i_1}=b_{i_1},x_{i_2}=b_{i_2},\ldots,x_{i_k}=b_{i_k}}$.\end{definition}
Thus, the integral of an $n$-form $\omega = f d_1 \wedge d_2\ldots \wedge d_n$ over the Hamming cube is obtained by summing the value of $(x_1 x_2\ldots x_nf)$ evaluated at all points on the Hamming cube. Trivially, it is seen that in fact $\int_{\mathbb{H}^n}\omega =f|_{x_1=1,x_2=1,\ldots,x_n=1}$. 
\begin{example}\label{ex:1}
For instance, let $n=3$, and consider the $1$-form $\omega=f(x_1,x_2,x_3) d_1.$ We have $\int_{supp(\omega)} \omega = \sum_{b_1,b_2,b_3 \in \mathbb{F}_2} x_1f|_{x_1=b_1,x_2=b_2,x_3=b_3} = \sum_{b_2,b_3 \in \mathbb{F}_2}f|_{x_1=1,x_2=b_2,x_3=b_3}.$
\end{example}
An $(n-1)$-form can be integrated over any $(i,j)$ boundary $\mathbb{H}^{n-1}_{i,j}$ of $\mathbb{H}^n$, with the tacit assumption that integrating over a face not contained in the support of the $(n-1)$-form yields a value of $0$. We can then define the integral of an $(n-1)$-form over the entire boundary $\partial \mathbb{H}^n$ of $\mathbb{H}^n$ as follows.
\begin{definition}Let $\omega$ be an $(n-1)$-form. Define $\int_{\partial \mathbb{H}^n}\omega =\sum_{i=1,\ldots,n, j\in \mathbb{F}_2}\int_{\mathbb{H}^{n-1}_{i,j}}\omega$.\end{definition}

\begin{example}\label{ex:2}
As an example, let $\omega=f d_1 \wedge d_2 \ldots \wedge \hat{d}_k \wedge \ldots d_n$ where, as before, $\hat{d}_{k}$ indicates the absence of the quantity $\hat{d}_{i}$. Evaluating the integral of $\omega$ over an $(i,j)$ boundary we see that $\int_{\mathbb{H}^{n-1}_{i,j}}\omega = 0 $ if $i\neq k$, since its support are the faces $\mathbb{H}^{n-1}_{k,0},\mathbb{H}^{n-1}_{k,1}$. If $i=k$, we have $\int_{\mathbb{H}^{n-1}_{i=k,j}}\omega = f|_{x_1=1,x_2=1,\ldots,x_k=j,\ldots,x_n=1}$. Therefore, we see that $\int_{\partial \mathbb{H}^n}\omega =\sum_{i=1,\ldots,n, j\in \mathbb{F}_2}\int_{\mathbb{H}^{n-1}_{i,j}}\omega = \int_{\mathbb{H}^{n-1}_{k,0}}\omega+\int_{\mathbb{H}^{n-1}_{k,1}}\omega$.\end{example}
We come now to the main theorem of this paper, a Stokes theorem like result: the Stokes-Zhegalkin theorem.

\begin{theorem}\label{thm:SZ}Let $\omega$ be an $(n-1)$-form. Then, 
\begin{equation}\int_{\mathbb{H}^n}d\omega=\int_{\partial \mathbb{H}^{n}}\omega.\end{equation}
\begin{proof}
We prove the theorem for an $(n-1)$-form of the kind $\omega=f d_1 \wedge d_2 \ldots \wedge \hat{d}_k \wedge \ldots d_n$. Given $\omega=f d_1 \wedge d_2 \ldots \wedge \hat{d}_k \wedge \ldots d_n$, we compute $d\omega=df \wedge d_1 \wedge d_2 \ldots \wedge \hat{d}_k \wedge \ldots d_n$. Now, $df = \sum_{i=1}^{n}\partial_i(f) d_i$, and therefore, $$d\omega = (\sum_{i=1}^{n}\partial_i(f) d_i) \wedge d_2 \ldots \wedge \hat{d}_k \wedge \ldots d_n = \partial_k(f) d_1 \wedge d_2 \wedge \ldots d_n.$$
By proposition \ref{prop:I1}, we have $\partial_k(f) = f|_{x_k=0}+f|_{x_k=1}$, so we conclude $d\omega = (f|_{x_k=0}+f|_{x_k=1})d_1 \wedge d_2 \wedge \ldots d_n$. Now, we have $$\int_{\mathbb{H}^n}d\omega = (f|_{x_1=1,x_2=1,\ldots,x_k=0,x_{k+1}=1,\ldots,x_n=1}+f|_{x_1=1,x_2=1,\ldots,x_k=1,x_{k+1}=1,\ldots,x_n=1}).$$
We also have from example \ref{ex:2} that $\int_{\partial \mathbb{H}^n}\omega = \int_{\mathbb{H}^{n-1}_{k,0}}\omega+\int_{\mathbb{H}^{n-1}_{k,1}}\omega = f|_{x_1=1,x_2=1,\ldots,x_k=0,x_{k+1}=1,\ldots,x_n=1} + f|_{x_1=1,x_2=1,\ldots,x_k=1,x_{k+1}=1,\ldots,x_n=1}$. We therefore conclude that $\int_{\mathbb{H}^n}d\omega=\int_{\partial \mathbb{H}^{n}}\omega$.
\end{proof}
\end{theorem}

\section{Summary}
This paper provided a unified theory of Boolean calculus culminating in the Stokes-Zhegalkin theorem, and thereby offering a complete perspective of the works \cite{encinas,tucker,StPosBook1}. It is possible to consider some extensions of the theory presented here. For instance, one can now consider the de Rham complex of a Boolean algebra and glean valuable insight from the associated cohomology groups. Also, the theory of Boolean differential equations \cite{StPosBook2,sefarati95} can be re-interpreted from the general framework presented here.

\bibliography{zbcalcref}
\bibliographystyle{plain}

\end{document}